\documentclass{article}
\usepackage{graphicx}
\usepackage{cite}
\usepackage{amssymb,amsmath}
\usepackage{layout}
\usepackage{subcaption}
\usepackage{amsthm}
\usepackage{algorithm}
\usepackage{algpseudocode}
\usepackage[makeroom]{cancel}
\usepackage[dvipsnames]{xcolor}

\DeclareGraphicsExtensions{.eps,.jpg,.pdf,.ps} 

\graphicspath{{./Figs/}}

\def\printfig{1}   % Set to "0" to skip figures, "1" to print them
\def\appendixa{1}  % Set to "0" to skip appendix a, "1" to include 
\def\appendixb{1}  % Set to "0" to skip appendix b, "1" to include
\def\appendixc{0}  % Set to "0" to skip appendix c, "1" to include - not yet used

\begin{document}

\title{Convergence analysis of Left-Right splitting surface scattering
  method}  
\author{Paul E Barbone{\footnote{Department of Mechanical Engineering,
      Boston University,  MA 02215, USA}},~  
{Mark Spivack\footnote{Centre for Mathematical Sciences, 
   University of Cambridge, CB3 0WA, UK}} ~
and Orsola Rath Spivack$^\dagger$
}  
 
\font\eightrm=cmr8

\date{\today}

\maketitle

\centerline{ Email: \texttt{barbone@bu.edu},
  \texttt{ms100@cam.ac.uk}, and
  \texttt{or100@cam.ac.uk}}

\medbreak

\begin{abstract}
We study the convergence of the Left-Right
splitting method (equivalent in key respects to the Method
of Multiple Ordered Interactions and Forward-Backward method) for wave
scattering by rough 
surfaces.  This is an operator series method primarily
designed for low grazing incidence and found in many cases to
converge rapidly, often within one or two terms even for large incident angles.
However, convergence is not guaranteed and semi-convergence may be observed.

Our aims are two-fold: (1) To obtain theoretical and physical insight
into the regimes in which rapid convergence occurs and the mechanisms
by which it fails, by examining and modifying eigenvalues of the
operator; (2) provide a strategy for increasing the speed of
convergence or more crucially for overcoming divergence, and providing
a stopping criterion. The first is addressed by subtracting successive
dominant eigenvectors from the incident field, to examine the impact
on divergence and on the incident spectrum.  For the second, we apply
a generalisation of Shanks' transformation to the operator series;
this effectively improves convergence and (unlike eigenvalue
subtraction) readily generalises to 3D and composite problems. These
results also explain why the method converges so rapidly for much
larger incident angles.  Finally we ask and give an analytical
solution to a key question: For a divergent eigenvector of the
iterating operator, what is the exact solution and can it be deduced
from the divergent series?  We show that the exact solution is
well-behaved and can be found from the series in a way which is
related to the Shanks transformation.

\end{abstract}

\def\EQN{\eqno}
\def\beq{\begin{equation}}
\def\eeq{\end{equation}}
\def\beqn{\begin{eqnarray}}
\def\eeqn{\end{eqnarray}}
\def\bef{\begin{figure}}
\def\eef{\end{figure}}
\def\del{\partial}
\def\Psiz{\Phi}
\def\psinc{\psi_{inc}}
\def\vnc{\v_{inc}}
\def\r{{\bf r}}
\def\rr{{\bf r'}}
\def\d{{\bf d}}
\def\h{{\bf h}}
\def\A{{\bf A}}
\def\B{{B}}   %{{\bf B}}
\def\bea{\begin{aligned}}
\def\eea{\end{aligned}}
\def\zz{Z}
\def\Arcsin{{\rm Arcsin}}
\def\Arctan{{\rm Arctan}}
\def\t{\theta_n}
\def\tt{\frac{\theta_n}{2}}
\def\hmax{h_{\max}}
\def\hmin{h_{\min}}
\def\P{{P}}   %{{\bf P}}
\def\Q{{Q}}     %{{\bf Q}}
\def\Fcal{\cal{F}}
\def\Gcal{\cal{G}}
\def\Scal{\cal{S}}
\def\Einc{{E_{inc}}}
\def\Escat{{E_s}}
\def\Hinc{{H_{inc}}}
\def\Hscat{{H_s}}
\def\Epartial{{\partial E\over \partial n}}

\def\Sn{{S_n}}
\def\Sone{{S_{n+1}}}
\def\Stwo{{S_{n+2}}}
\def\v{{\bf v}}
\def\u{{\bf u}}
\def\w{{\bf w}}

\def\rbar{{\bf r}}
\def\lbar{{\bar{\lambda}}}
\def\J{{\bf J}}
\def\t{{\bf t}}

\parindent 0 true cm

\section{Introduction}\label{intro}

\smallbreak

We study the convergence of the Left-Right
splitting method (effectively equivalent to Method of Multiple
Ordered Interactions; Forward-Backward method) for wave scattering by rough
surfaces
\cite{Kapp,pino,colak2,tran1997calculation,spivack2017efficient,spivack2001validation}.
 
This is an operator series method primarily
designed for low grazing incidence and found in many cases to
converge rapidly, so that just one or two terms are needed in many
cases even for larger incident angles. 
However, convergence in general is not guaranteed, and in some cases an
apparently converged solution gives way to a rapidly growing
nonphysical component. 

Much of this analysis is for 2D problems but the generalisation to
3D is straightforward. 

Our aims are two-fold: (1) To obtain theoretical and physical insight
into the regimes in which convergence is slow and mechanisms by which
it fails, by examining eigenvectors of the iterating operator and
deriving the exact solution, and by examining and solving exactly
for diverging eigenvectors of the iterating operator;
(2) provide a strategy for overcoming
divergence or improving convergence.  By subtracting successive
dominant eigenvectors from the incident field we can slow or reverse
the divergence and examine the incident spectrum.  Further, we apply a
generalisation of Shanks' transformation to the operator series; this
is more effective, and (unlike the first approach) readily generalises
to 3D and composite problems.

We focus mainly on the case of TM (vertically polarised)
electromagnetic incident field; equivalent to Neumann boundary
condition for acoustic waves. By comparison with the TE (Dirichlet)
case this gives rise to greater diagonal dominance of the discretised
operator. and therefore more readily generalises to 3-dimensional 
problems.  Convergence is remarkably robust with respect to incident
angle, but unsurprisingly can be overwhelmed by sufficiently rough
surfaces, i.e. large spatial surface variance.  A critical issue to
address is optimal \textit{stopping point} since some cases give good
results at the first term but subsequently diverge, i.e. they exhibit
semiconvergence. 

In summary, the key findings are: (a) Convergence can occur for any
incident angle for moderately rough surfaces, while semi-convergence
often occurs even for very rough surfaces; (b) as might be expected,
divergence of the series is typically due to dilating (dominant)
eigenvectors of the iterative operator, and worsens with greater
surface variation; (c) for a divergent incident eigenvector, the full
integral equation is well-behaved and we can obtain its solution
exactly; (d) the incident data residual $||A\psi_n - \psinc||$ is
well-behaved and can therefore be used as a stopping criterion; (e) we
formulate scalar and vector versions of Shanks' transformation and
find they can be used to improve convergence at low computational
cost.  We also examine the relationship between eigenvector
subtraction and Shanks iterations.

This paper is organised as follows: Section \ref{math}: mathematical
formulation and operator series solution; 
Section \ref{results}: basic numerical results; Section
\ref{eigvecresults}: projection methods; Section \ref{shanks}: scalar
and vector Shanks' transformation.

\begin{figure} 
\hskip 2.8 true cm    
\if\printfig1\includegraphics[width=0.6\linewidth]{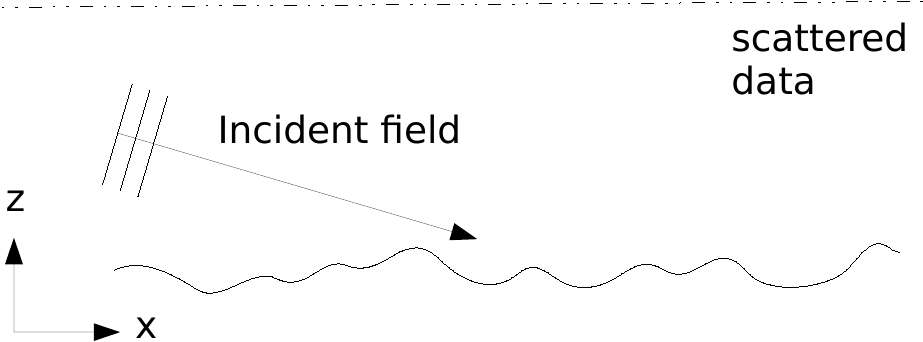} \fi
    \caption{Schematic view of scattering configuration} 
  \label{schematic} 
\end{figure}

\section{Mathematical formulation and L-R splitting
  series} \label{math} 

Coordinate axes $x$ and $z$ are taken as in Fig. \ref{schematic}
where $x$ is the horizontal and $z$ is the vertical.
We consider scattering of a time-harmonic wavefield incident from the
left on an extended perfectly reflecting rough surface $h(x)$ 
with mean plane lying along $x=0$.
We will assume the field is incident at low
grazing angles resulting in small angles of scatter.  
The total field $E$ obeys the Helmholtz equation 

\begin{equation}
\nabla^2 E(x,z)+k^2E(x,z)= 0 %f(x,z)
\end{equation}
%{
where $k$ is the wavenumber.
 The  2-dimensional free space Green's
  function $G_0$ is the zero order Hankel function of the first kind,
  $G_0(\r,\r')=(1/4i)H_0^{(1)}(k|\r-\r'|)$.%}

The field in the medium can be written as 
a boundary integral over the normal derivative along the surface.
The incident electric field is assumed to be time-harmonic, with 
time-dependence $\exp(-i\omega t)$, say, and
may be taken to be horizontally or vertically plane
polarized, i.e. corresponding to TE or TM.  
We now suppress the time-dependence and consider the time-reduced
component, and will initially assume an incident TM field.
Suppose that the wave $H$ is scattered by a rough perfectly conducting 
one-dimensional surface $h(x)$,
so that $H$ obeys the Helmholtz 
wave equation $(\nabla^2+k^2)H=0$.  This is shown schematically in
figure 1.  
Let $G$ be the free space 
Green's function, so that $G$ is the zero order Hankel function of the
first kind,
\beq 
G(\rbar,\rbar') = {1\over 4i} H_0^{(1)} (k|\rbar-\rbar'|) .
\eeq
In this case the normal derivative $H_n$ of the field at the surface 
vanishes, and the governing integral equation is then obtained as

\beq
\Hinc(\rbar_s) = H(\rbar_s) - \int_{-\infty}^\infty 
{\del G(\rbar_s,\rbar')\over \del n}
H (\rbar') dS
\label{(3.1)}
\eeq
where $n$ denotes the outward (i.e. downward) normal, 
integration is over the
surface, and
$\rbar_s=(x,h(x))$, $\rbar'=(x',h(x'))$ both lie on the surface.
(The integral in (\ref{(3.1)}) must be interpreted with care, since 
$\del G/\del n$ is singular at $\rbar'=\rbar_s$. This expression corresponds 
to the limit of the boundary integral for a point $\rbar$ tending to the 
surface, $\rbar \rightarrow \rbar_s$.)
For convenience we write eq.(\ref{(3.1)}) in operator notation,

\beq
\Hinc (\rbar_s) = (L+R)H \label{(3.2)}
\eeq
with corresponding field integral

\beq
\Hscat (x,z) = -(L+R)H \label{(3.3)}
\eeq
where $L$ and $R$ are defined by

\beq
Lf(x,z)=- \hspace{-11pt}
\int_{-\infty}^{x} {\del G(\rbar_s,\rbar')\over \del n} f(x') dS,
~~~~~ Rf(x,z)=\int_{x}^\infty  {\del G(\rbar_s,\rbar')\over \del n} f(x') dS
\label{(3.4)}
\eeq
and $\rbar=(x,z)$, $\rbar'=(x',h(x'))$, and $L$ includes the principal
value of the integral.

Integral equation (\ref{(3.2)}) has formal solution

\beq
H  = (L+R)^{-1}\Hinc   \label{(3.5)}
\eeq
which can be expanded in a series

\beq
H
= \left[L^{-1} - L^{-1} R L^{-1} +L^{-1}  \left(R L^{-1} \right)^2 -
...+  (-1)^n L^{-1}\left(R L^{-1} \right)^n ...   \right] \Hinc . \label{(3.6)}
\eeq

We will denote the operator $\B = - R L^{-1}$  so that (\ref{(3.6)})
can be written   
\beq
H
~~ = ~~ L^{-1} \sum_{n=0}^\infty B^n \Hinc ~~
. \label{(3.6b)}
\eeq

Provided it converges this series can be truncated, and treated term by term.
When the system is discretized, the operator $L$ yields
a lower triangular matrix. Similarly $R$ becomes upper triangular (with zero 
on the diagonal).
Inversion of the matrix $L$ can be carried out very 
efficiently (using Gaussian elimination and backward substitution)
to give the first term of eq.(\ref{(3.6)}).
Since subsequent terms in the series are products of $L^{-1}$ and $R$, they
can also be evaluated efficiently. 

The heuristic argument for the convergence of the series is as follows:
Convergence of the series can be expected if the effect of the operator
$R$ on its argument is `sufficiently small' compared with that of $L$
(although, as with most scattering approximations, it is
difficult to put rigorous bounds on the surface statistics which ensure 
convergence).   Note that as surface roughness increases $R$ itself may no
longer be considered small since its norm may become comparable with
that of $L$.  However, for predominantly right-going waves, the
functions on which $R$ operates in the series (\ref{(3.6)}) will have
phase components varying 
rapidly with $x$,  for which the right half-integral 
represented by $R$ will give rise to functions whose
amplitude is small, as required.  
The purpose of the present study is in part to examine this argument
in more detail and consider regimes  
where it breaks down.

\section{Numerical results}\label{results}

We begin by giving numerical examples, in order to motivate the
investigation of convergence which follows.  Rough surfaces are
generated as in \cite{spivack2001validation} in which the solutions
were carefully validated against finite element time domain BAE
Systems models, scaled as necessary here to induce divergence.  For
this purpose, the algorithm was applied to rough surface patches
embedded on longer and otherwise flat surfaces as in
\cite{spivack2001validation}, and also to surfaces which were rough
throughout their extent.

The key computational advantage of the L-R method is that the
inversion of a full matrix is replaced by a small number of inversions
of lower-triangular matrix $L$,  which can be solved by Gaussian
elimination. (In 3-dimensional problems $L$ becomes lower
block-triangular but the principle is similar.)

The main computational expense in 2D thus arises from solving systems
of equations of the type 
$L\v=\v_{inc}$) achieved without an explicit matrix inversion or
eigenvalue calculations. 
(The secondary expense is `matrix filling' to evaluate the entries of
$L$ and $R$.  
In 3D the computational cost of these tasks is in practice reversed.)
Here however 
we are interested explicitly in such properties of the matrices, and
for simplicity the bulk of the calculations are carried out in MATLAB.  

More explicitly: In 2D for $n$ range steps, without further
optimisations, using k steps of the series (\ref{(3.6b)}),  inversion
accounts for $O(kn^2)$ operations and matrix filling
for $O(n^2)$. 
Here $k$ is typically small (perhaps 1 or 2).    [In 3D with $n$
range steps and $m$ transverse points, inversion takes around $O(kn^2
m^3))$ at most whereas matrix evaluations around require
$O(n^2m^2)$. At first glance this suggests inversion should again
dominate, but in practice the matrix evaluations have a high multiple,
and the $m^3$ term can be reduced to less than $m^2$.]

\begin{figure} [ht]
	\hskip 2.8 true cm    
\if\printfig1	\includegraphics[width=0.6\linewidth]{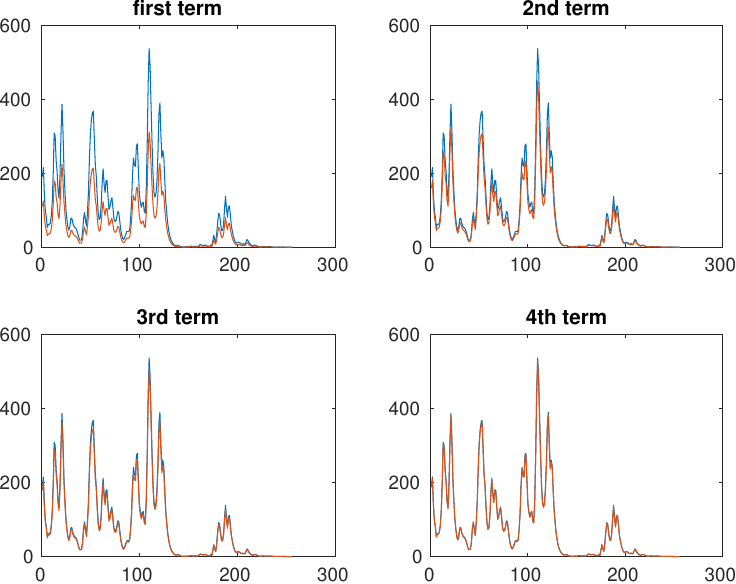} \fi
	\caption{Example of convergence of L-R method for typical regime} 
	\label{fig2} 
\end{figure}

Two important and related observations here:  \\
(1) Surface variation and its effect on the operator appears far more
significant than incident angle, for convergence properties.\\  
(2) In well-behaved cases the series converges even for very large
incident angles (near-normal or even backward-going). \\ 
%{\bf
  We will show below that for a range of surface roughness parameters
convergence may be  
observed for any incident angle. %}

\section{Divergent eigenvectors and eigenvector projection}\label{eigvecresults}

In order to investigate divergence and the regimes in which it occurs
[{see section \ref{results}}] we increased surface
heights and incident angles to produce divergence in the series
solution (\ref{(3.6)}). When the series diverges a consistent feature
was observed: for large $n$ the $n$-th term of the series becomes
dominated by (say) $\lambda^n\v_0$ where $\lambda$ is some scalar with
$|\lambda| > 1$ and $\v_0$ is some vector.  This suggests that
$(L\v_0,\lambda)$ is an eigenvector/eigenvalue pair for the operator
$\B$, and both can be found approximately from the numerical solution.
We will refer to such eigenvectors as {\sl
  dilating}.

This raises several questions. In particular: (1) Can we determine the
exact solution due to a given dilating vector $\v_0$ for the series
solution?  (2) Can we subtract the divergent contribution due to
$\v_0$, or eliminate it from the data to produce a better-conditioned
series; and if so can we extend the principle to multiple dilating
eigenvectors?

\textit{Remark:} Since the eigenvectors of
self-adjoint operators form an orthonormal set, they can be efficiently
identified and subtracted one-by-one, for example using the power
method.  However, the operator $B$ here is non-self-adjoint; thus
eigenvector subtraction requires calculation of the same computational
order as full inversion, and therefore not useful in practice.  (See
below and Appendix A for summary.)

\subsection{Projection method}\label{projectionmethod}

One goal here is to investigate the role of dilating eigenvectors in
determining the convergence properties of the series.  To this end we
aim initially to identify such components from the initial data,
calculate the\textit{ convergent contribution} due to these as
explained below, and then recalculate the series. We can verify this
approach by directly examining the full solution applied to these
components.

\bigbreak

Suppose we have a complete set of eigenvectors $\v_i$ for $B = -RL^{-1}$ with 
eigenvalues $\lambda_i$.  We assume for convenience that the eigenvectors have 
multiplicity one.   These eigenvectors do not form an orthogonal set 
(as would be the case if $B$ were self-adjoint).

Denote by $P_i$ the rank 1 projection  onto the eigenspace of $\v_i$, so that
$P_iP_j = 0$ for $i\neq j$ and 
the sum of any subset of $\{P_i\}$ is a projection. However the range and null space of 
any $P_i$ are not orthogonal subspaces. 

For a general vector $\v$ consider the partial sums of [finite approximations to] the 
operator series on the 
right-hand-side of equation (\ref{(3.6b)}), say:
\beq S_n (\v)=  L^{-1}\sum_{k=0}^n B^k \v \eeq 
This can be written in terms of the projections $P_i$ as

\beq S_n (\v)~=~ L^{-1}\sum_i \sum_{k=0}^n \left(\lambda_i^k P_i \right)\v 
~=~  L^{-1}\sum_i\sum_{k=0}^n \left(\lambda_i^k \alpha_i\v_i \right)\eeq 
where $P_i\v = \alpha_i\v_i$, which gives
\beq 
S_n (\v)~=~ L^{-1}\sum_i \alpha_i \Big(\frac{ \lambda_i^{n+1} - 1
}{\lambda_i - 1} \Big)~\v_i   . 
\label{snv}
\eeq

Now apply this to a single eigenvector $\v_i$, in order to compare with the 
exact solution. Provided $\lambda_i \neq 1$, equation (\ref{snv}) reduces to   
\beq
L^{-1} \sum_{k=0}^n B^k \v_i
~=~ L^{-1} \sum_{k=0}^n \lambda_i^k \v_i
~=~  \left(\frac{\lambda_i^{n+1}-1}{\lambda_i-1}\right)~L^{-1}~\v_i
\label{finitesum}
\eeq
which clearly converges for $|\lambda_i| < 1$ and diverges for $|\lambda_i|>1$.

Now consider the exact solution for $\v_i$:
\beq (1-\lambda_i)\v_i ~~=~~\v_i + RL^{-1}\v_i ~~=~~ (L+R)L^{-1}\v_i
~~=~~ AL^{-1} \v_i 
\label{lemma1a}\eeq 
Thus the integral equation (\ref{(3.5)}) is well-behaved at $\v_i$,
with exact solution given  by 
\beq A^{-1}\v_i ~~=~~ \frac{1}{1-\lambda_i}L^{-1}\v_i  . \label{lemma1b}\eeq

Comparing with (\ref{finitesum}) we see that for $|\lambda_i| < 1$ the
series converges to the expected value,
but in either case the correct limit is obtained by
discarding the term $\lambda_i^{n+1}$ in the numerator. This result
extends to any sum of eigenvectors.  (Note also that it extends without  
modification to the analogous 3-dimensional scattering problem.) So for any $\v$ 
\beq A^{-1}\v ~=~
\sum_i \Big(\frac{ 1}{1-\lambda_i} \Big)~L^{-1} P_i \v_i. 
\label{exact}\eeq

\vskip 1 true cm

\begin{figure}[ht]
	\hskip 2.8 true cm    
\if\printfig1
\includegraphics[width=0.6\linewidth]{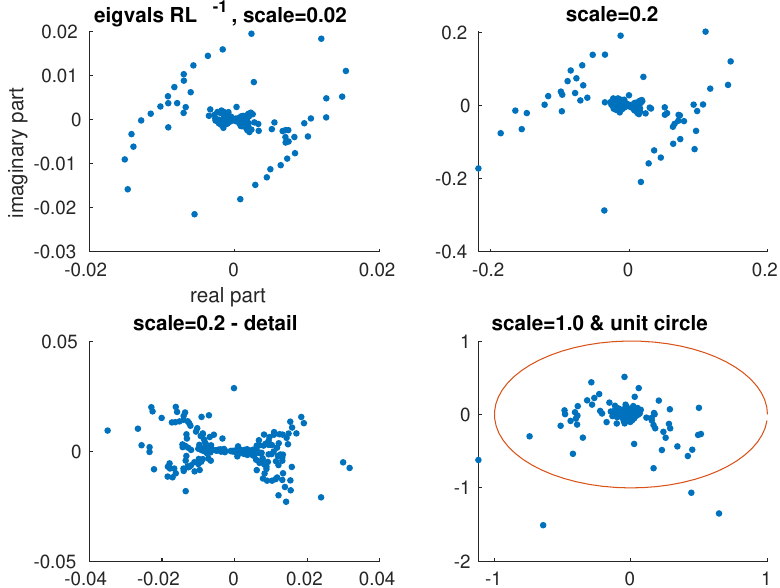} \fi 
	\caption{Eigenvalue distributions for increasing surface
          roughness; unit circle shown in largest-roughness example}  
	\label{fig3} 
\end{figure}

\begin{figure} [ht]
	\hskip 2.2 true cm    
\if\printfig1
\includegraphics[width=0.6\linewidth]{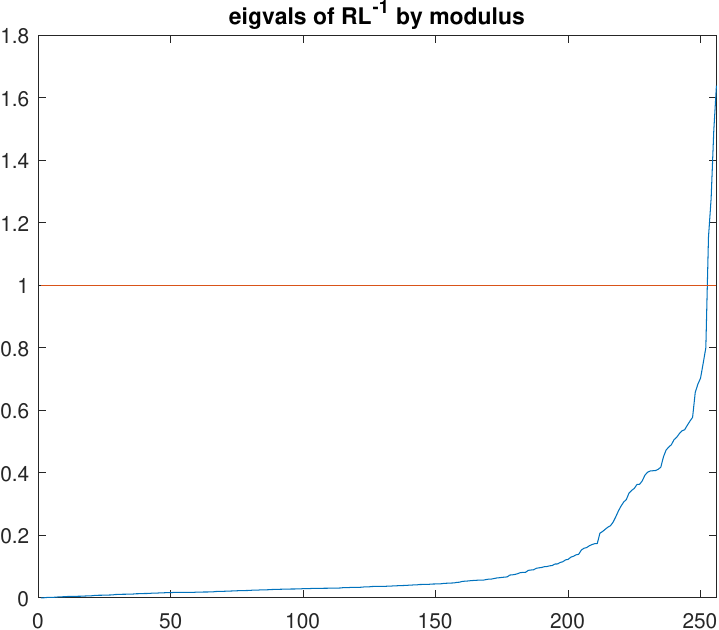} 	\fi  
	\caption{modulus $|\lambda|$ of eigenvalue of operator $B$ for
          'typical' divergent case}  
\label{fig3b} 
\end{figure}

\begin{figure} [ht]
	\hskip 2.2 true cm    
\if\printfig1	\includegraphics[width=0.6\linewidth]{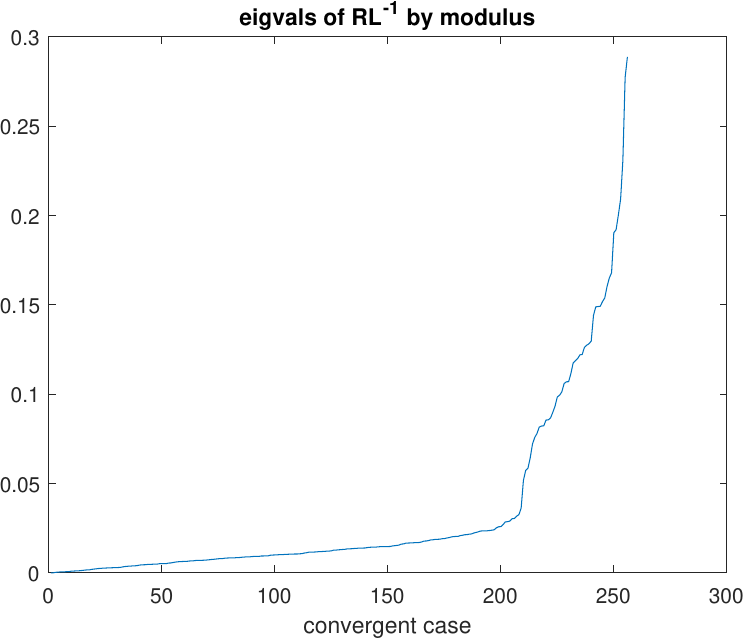} \fi
	\caption{modulus $|\lambda|$ of eigenvalue of operator $B$ for
          typical convergent case}  
	\label{fig3c} 
\end{figure}

\begin{figure} [ht]
	\hskip 2.0 true cm    
\if\printfig1
\includegraphics[width=0.6\linewidth]{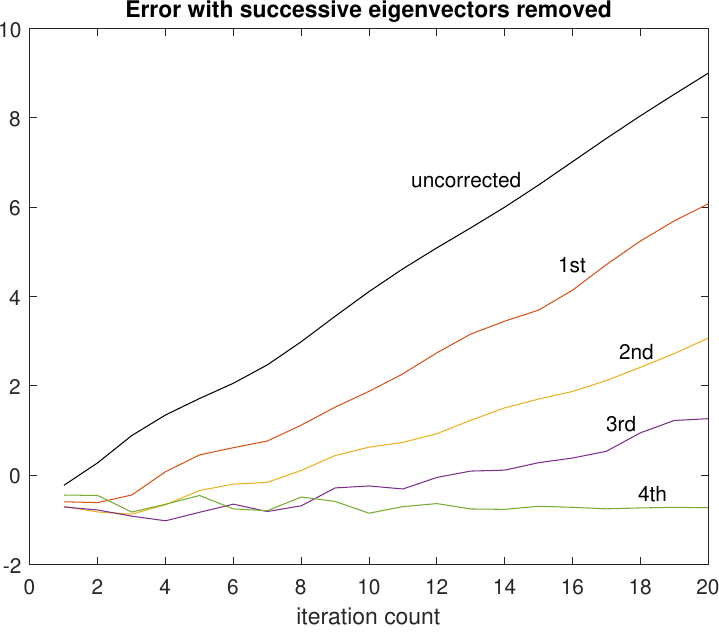}  
\fi
	\caption{Comparison between error in uncorrected L-R iteration
          and after removal of successive dominant eigenvectors,
          i.e. those whose eigenvalues have modulus greater than 1}  
	\label{eigvecsub} 
\end{figure}

\section{Shanks' transformation}\label{shanks}

 \subsection{Motivation} 

The Shanks transformation \cite{shanks,brezinski,bender} is a
remarkably effective method for accelerating the convergence of
slowly-convergent series. Although applied originally to scalar series
several versions have been proposed for vector and matrix series.   It
works roughly as follows:   

 Let S denote the limit of a series of
 partial sums, $S_n$.  Suppose a series behaves as $S_n \sim S +
 r^n$.  As long as  
 $|r|< 1$, then  $S_n \rightarrow S$ as $n \rightarrow\infty$.  
Regardless of $|r|$, however, $S_n - S_{n-1} = r^n - r^{n-1}$.  One may
thus solve for $r$, subtract its contribution to the sum, and get an improved
estimate of $S$.

\subsubsection*{Scalar Shanks transformation}

First suppose we have a sequence \Scal=$\{S_n\}$ which in most cases
will arise from the partial sums of a sequence (say) 
$\{a_n\}$.  

Suppose that
 \beq \Sn = S + \alpha \rho^n  \label{shanks1} \eeq
\beq {\Sone} = S +\alpha \rho^{n+1}  \label{shanks2} \eeq
\beq \Stwo = S + \alpha \rho^{n+2}  \label{shanks3} \eeq

From equations (\ref{shanks1})-(\ref{shanks3}) we can define:
\beq \Delta \Sn = \Sone -\Sn = \alpha(\rho^{n+1}-\rho^n) \eeq
\beq \Delta \Sone = \Stwo -\Sone = \alpha(\rho^{n+2}-\rho^{n+1}) \eeq

Then \beq \rho =  \frac{\Delta \Sone}{\Delta \Sn}, ~~~~ 
\alpha = \frac{\Delta\Sn}{\rho^{n+1}-\rho^n}   \label{scalarparams}\eeq
and we can simplify to obtain

\beq  S = \frac{\Sn \Stwo - S^2_{n+1}}{\Stwo - 2\Sone + \Sn} . \eeq

In general if equations (\ref{shanks1})-(\ref{shanks3}) are not exact
but valid asymptotically, this gives rise to Shanks' transformation:
${\Fcal}(\Scal) =$ $\{F_n\}$ where 
\beq  F_n = \frac{\Sn \Stwo - S^2_{n+1}}{\Stwo - 2\Sone +
  \Sn}  \label{scalarshanks}\eeq 
so that $F_n$ is a new sequence.  
If in equations (\ref{shanks1})-(\ref{shanks3}) $|\rho| < 1$ then
$F_n$ may converge to $S$ faster than the 
original sequence $S_n$, often dramatically so.    If $|\rho| > 1$ and
we replace $S$ in  
(\ref{shanks1})-(\ref{shanks3}) by a convergent sequence whose limit
is $S$, then  
(\ref{scalarparams}) remains valid asymptotically.

\begin{figure} 
	\hskip 2.0 true cm    
\if\printfig1 \includegraphics[width=0.6\linewidth]{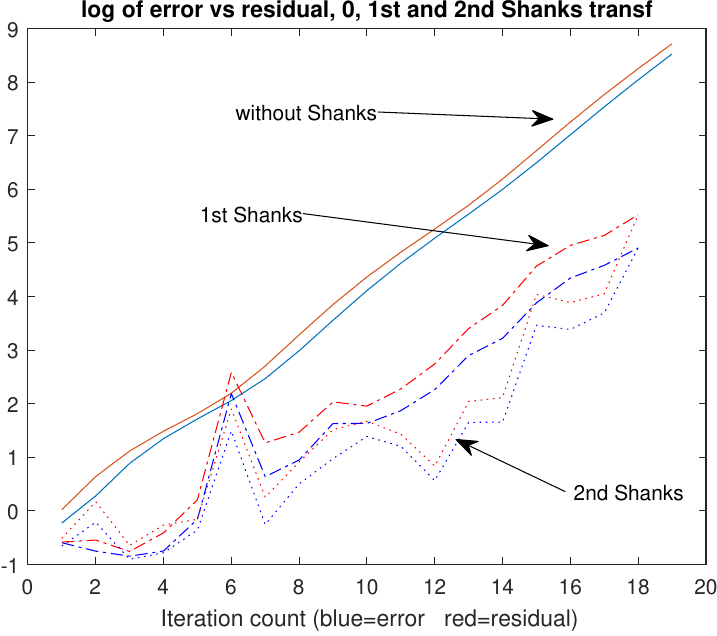} \fi
	\caption{This shows residual is good indication of error; Shanks improves solution 
	initially and typically best at 3rd iteration} 
	\label{fig5} 
\end{figure}

\begin{figure} 
	\hskip 2.0 true cm    
\if\printfig1 \includegraphics[width=0.6\linewidth]{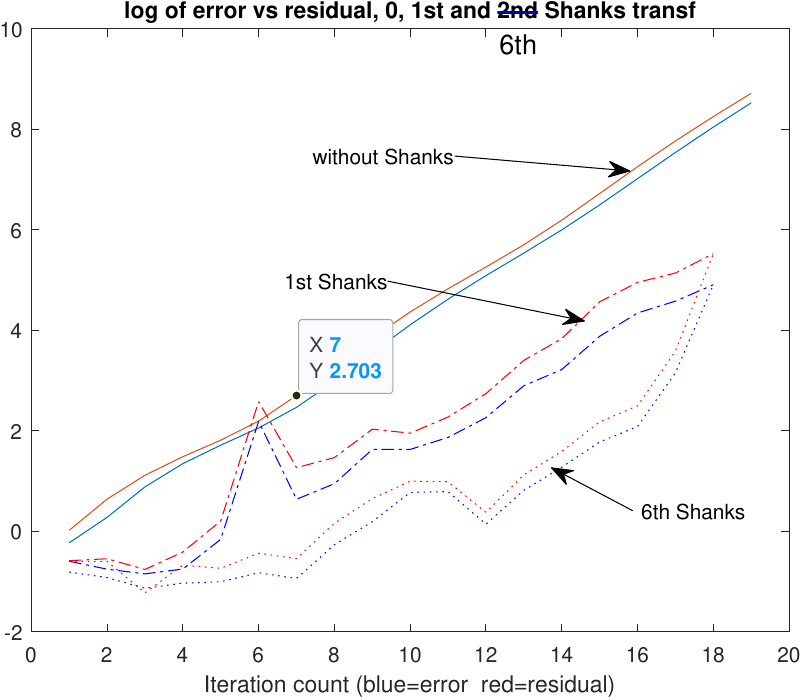} \fi
	\caption{Similar to Fig \ref{fig5}, for higher-order $
          $transformation. 6th Shanks is a
          significant further improvement and the repeated process has
          removed much of the variability [but see next figure].} 
	\label{fig6} 
\end{figure}

\begin{figure} 
	\hskip 2.8 true cm    
\if\printfig1
\includegraphics[width=0.6\linewidth]{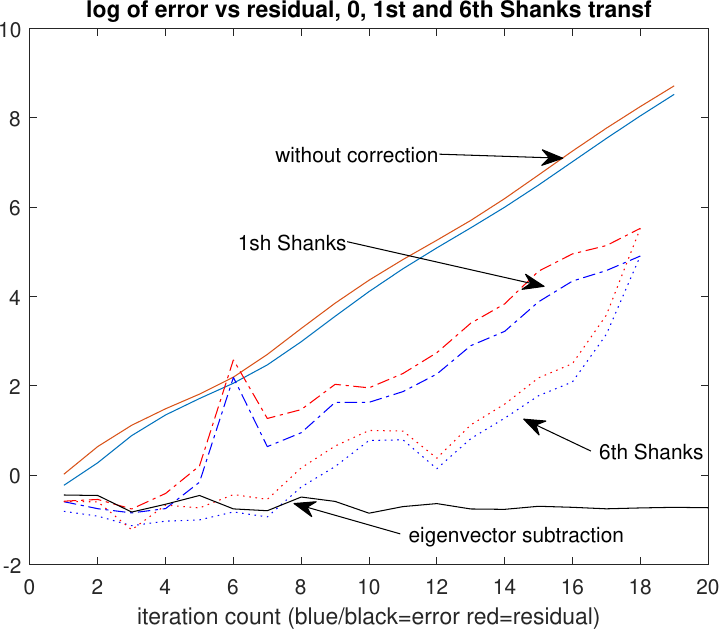}
\fi 
	\caption{Same as Fig \ref{fig6} with the error from
          eigenvector subtraction included.  Although eigenvalue
          subtraction is clearly more successful in preventing
          divergence, at its most accurate it does not  
   outperform Shanks transformation at low iteration count.}
	\label{fig6b} 
\end{figure}

\begin{figure} 
	\hskip 2.8 true cm    
\if\printfig1 \includegraphics[width=0.6\linewidth]{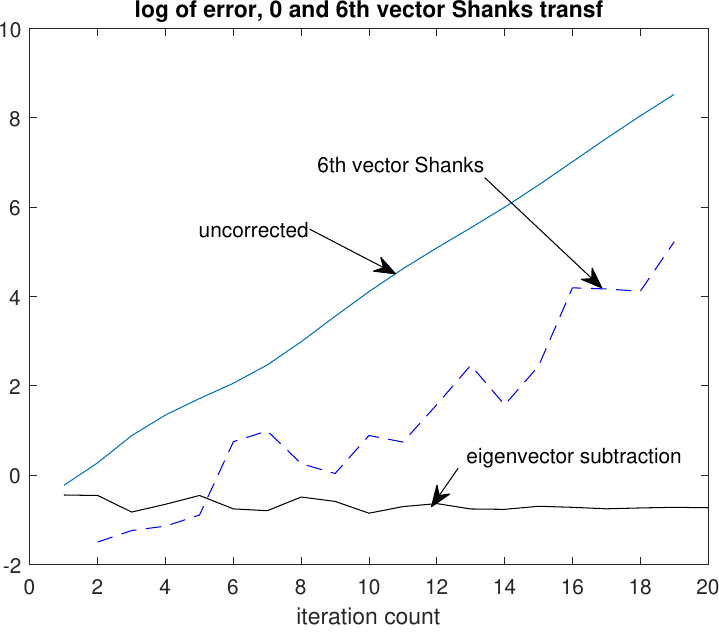} \fi
	\caption{Similar to Fig \ref{fig6}, using the vector form of
          the Shanks transformation} 
	\label{fig8} 
\end{figure}

\medbreak
\subsubsection*{Vector Shanks transformation}\label{vshanks}
Suppose now that $\Sn$, $S$ and $\v$ are vectors,  and that
 \beq \Sn = S +\lambda^n \v  \label{v1}\eeq
\beq {\Sone} = S +\lambda^{n+1} \v  \label{v2}\eeq
\beq \Stwo = S +\lambda^{n+2} \v  \label{v3}\eeq

Define:
\beq \Delta \Sn ~=~ \Sone -\Sn ~=~ (\lambda^{n+1}-\lambda^n)\v \label{v3b}\eeq
\beq \Delta \Sone ~=~ \Stwo-\Sone ~=~
(\lambda^{n+2}-\lambda^{n+1})\v ~=~ \lambda \Delta \Sn \eeq

Thus the vectors $\Delta \Sn$, $\Delta \Sone$ are colinear, so
that $\lambda$ can be obtained in terms of inner products, and  $\v$
directly from $\Delta \Sn$: 
\beq 
\lambda = \frac{ \left<\Delta\Sone, \Delta\Sn \right>}
     { ||\Delta\Sn||^2}. ~~~~~~ \v = \frac{ \Delta\Sn} {\lambda^{n+1}-\lambda^n}   .
\label{v4}     \eeq
We will not use $\v$ explicitly in this transformation, but we can
compare it - or rather its limiting value - with the dominant
eigenvector obtained in section \ref{eigvecresults}.

From $\lambda\times(\ref{v1})-(\ref{v2})$:
\beq \lambda S + \lambda^{n+1}\v - S -\lambda^{n+1}\v ~=~ \lambda\Sn - \Sone \eeq
\beq (\lambda-1)S ~=~\lambda\Sn - \Sone \eeq
which gives $S = {\lambda \Sn - \Sone}/(\lambda-1) $ with $\lambda$ given by (\ref{v4}).
When equations (\ref{v1})-(\ref{v3}) are approximate we allow
$\lambda$ to vary with $n$ to obtain vector Shanks' transformation:
${\Gcal}(\Scal) =$ $\{G_n\}$ where 
\beq 
G_n = \frac{\lambda_n \Sn - \Sone}{\lambda_n-1} \label{vectorshanks}
\eeq 
Notice that we may apply this repeatedly to obtain higher-order Shanks transformations ${\Gcal}^m$, for example
${\Gcal}^2(\Scal) = \Gcal(\Gcal(\Scal)$.

\medbreak
\subsection{Shanks transformations applied to L-R series:}

Now consider the L-R series (\ref{(3.6)}). $S_n$ now represents the
partial sums whose behaviour we wish to model. In this case for
sufficiently large $n$, and when the series starts to diverge we
observe behaviour in line with equations (\ref{v1})-(\ref{v3}).  It is
relatively easy to estimate $\lambda$ and $\v$  from the higher terms
of the numerical solution, but it is not clear how to remove them, for
reasons explained elsewhere in this paper.

We may first try to apply the scalar version of Shanks' transformation
(\ref{scalarshanks}) in a pointwise manner, treating the solution at
each point on the surface independently.

{This scalar approach has a couple of obvious flaws
  as it stands,  although it has been found to be 
  effective. (1) It makes no use of the spatial continuity or
  correlations of the underlying functions; (2) it implicitly solves
  for a different $\lambda$ at each surface point; and (3) the
  denominator in  (\ref{scalarshanks}) may occasionally be near zero
  and cause localised numerical errors.} 

To address some of these concerns we apply the vector form of Shanks
transformation (\ref{vectorshanks}).  At low iteration count the
result  improves upon the scalar version and  significantly upon the
uncorrected L-R series. Fig. \ref{fig8} shows the result of the 6th
order transformation in comparison with both the uncorrected L-R
series and the result of eigenvector subtraction.

Fig. \ref{fig8} shows the result of applying the vector form of the
Shanks transformation to the divergent case, compared with eigenvector
subtraction.  It is interesting to note that although it eventually
diverges, the 6th vector Shanks transform is significantly better at
small iterations, but the comparison is somewhat misleading, since for
each $n$ the first term $S_1^n$ of the $n$-th Shanks transformation requires
$n+1$  terms of the original series.
l

\begin{figure} 
	\hskip 2.8 true cm    
\if\printfig1
\includegraphics[width=0.6\linewidth]{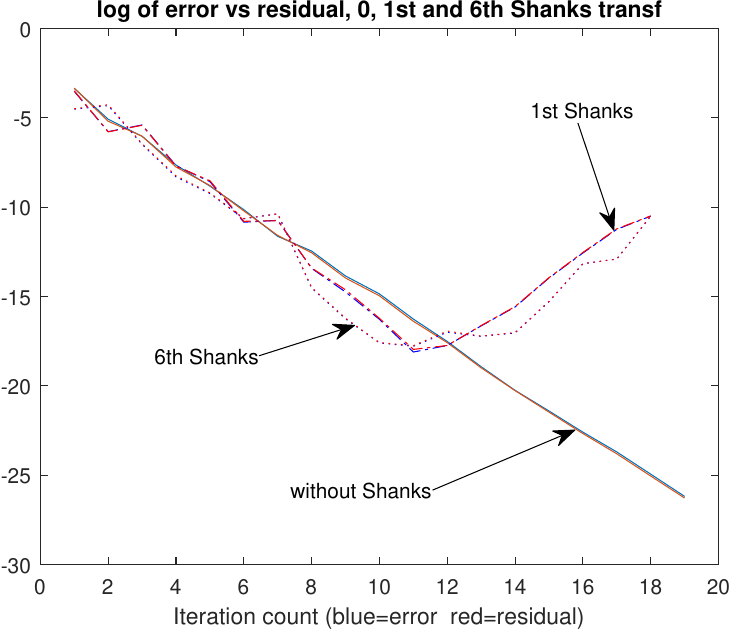}
\fi 
	\caption{Similar to Fig \ref{fig6}, applied now to a convergent case}
	\label{fig7} 
\end{figure}

Fig. \ref{fig7} shows the result of applying Shanks transformation
to a case in which convergence is rapid. The higher-order
transformation does appear to be more accurate the the original
calculation at the first term, but this is again somewhat misleading.

\subsection{Relationship between Shanks transformations and eigenvectors}

We examine here the correspondence between eigenvectors of $\B$ and
Shanks vectors which are implicit in the formulation of the Shanks
transformation (\ref{v4}).  Denote by $\v_1, \v_2, ...$ the limiting
vectors produced by successive applications of Shanks method, and by
$\w_i$ the eigenvectors of $\B$ in descending order of eigenvalue
modulus.

\begin{figure} 
	\hskip 2.8 true cm    
\if\printfig1 \includegraphics[width=0.6\linewidth]{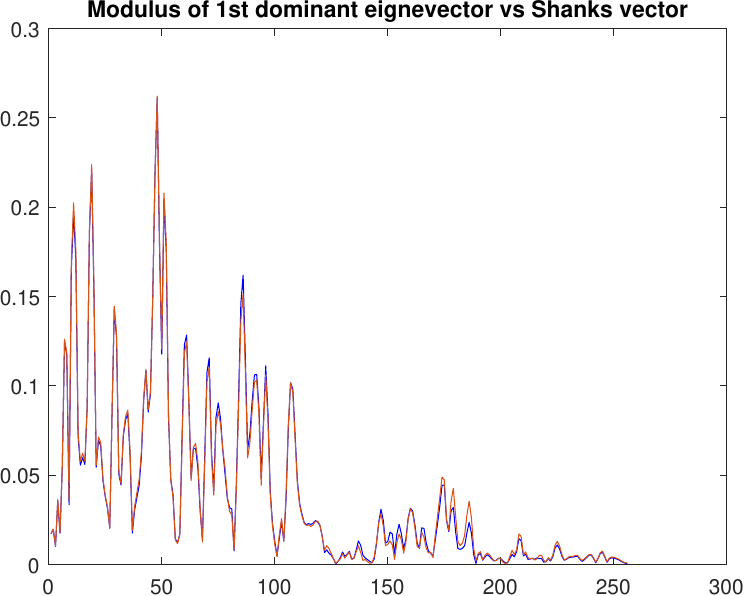} \fi 
	\caption{Eigenvector 1 vs Shanks$^1$ - close agreement}
	\label{eigvec_vs_shanks} 
\end{figure}

Once it has been normalised, the first order Shanks vector $\v_1$
closely reproduces the eigenvector $\w_1$, as shown in
Fig. \ref{eigvec_vs_shanks}; in other words they are colinear.
Although $\v_2\neq\w_2$, we find that $\w_2$ is coplanar with $\v_1,
\v_2$; and more generally the subspaces generated by the first $n$
vectors of each set are seen to coincide.

This explains why the Shanks method gives good results at low
iteration count, but eventually diverges in cases where the underlying
problem does so.
(See \cite{shanks,brezinski,bender,colak2007convergence,tran1997calculation})

\begin{figure} 
	\hskip 2.8 true cm    
\if\printfig1 \includegraphics[width=0.6\linewidth]{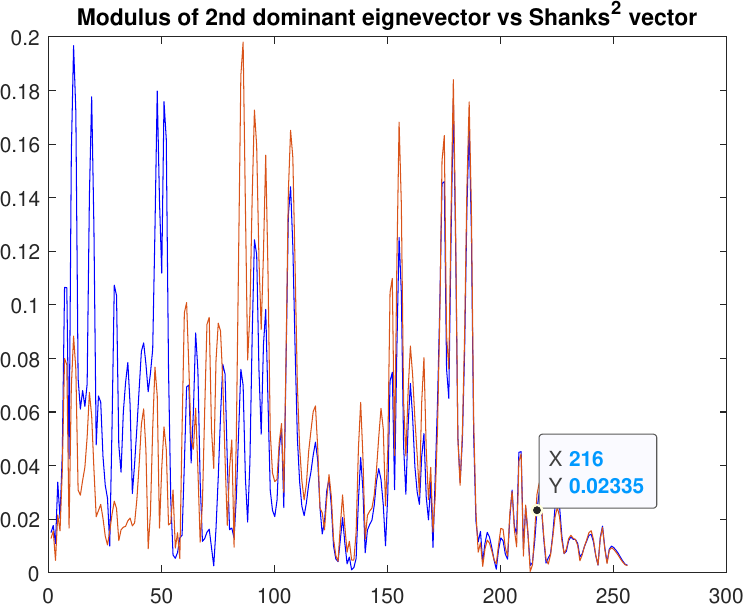} \fi
	\caption{Eigenvector 2 vs Shanks$^2$.}
	\label{eigvec2_vs_shanks2} 
\end{figure}

\section{Conclusions}

We have studied analytically and computationally the
convergence of the Left-Right splitting method for rough surface
scattering and the regimes in which it fails. We have proposed two
classes of method for overcoming divergence, based on eigenvector
subtraction and a form of Shanks' transformation.  Eigenvalue
subtraction illuminates the mechanisms by which convergence fails but
is not a practical method as it cannot feasibly be extended to
3-dimensional problems. On the other hand Shanks' transformation can
accelerate convergence and easily generalises to 3D.

The first few terms are typically well-behaved even when the method
ultimately diverges; in other words they display a type of
semi-convergence. This raised the question of finding a stopping
criterion.  We have shown that the residual (obtained by substituting
successive iterates into the forward scattering equation and comparing
with the initial data) closely tracks the solution error, and can thus
can be used as a reliable stopping criterion.

Finally for `dilating' or divergent eigenvectors of the iterating
operator $B$, we find that these are related to eigenvectors of the
exact scattering operator $A$. We can derive the eigenvalues exactly
from those of $B$, and find that we can estimate both eignevectors and
eigenvalues numerically from the divergent terms of the series.

In summary:
(1) We have analysed the convergence of L-R
splitting, and illuminated the regimes where it works well/badly.
This is necessarily carried out in representative example
cases.

(2) Eigenvector subtraction has been shown to succeed in removing
divergence (but is not usefully generalisable to 3D because
computational cost is of the same order as full inversion).

(3) L-R splitting has been shown to be applicable for very rough
surfaces and wide range of incident angles.   

(4) Scalar and vector Shanks transformation have been developed for the problem
and have been shown shown to accelerate convergence and in some cases (at low
iterations) remedy divergence.

(5) The issue of summing an optimum number of terms has been addressed
in semi-convergent cases  where the series initially converges and then
starts to diverge. We used residual 
$||A\psi_n - \psinc||$ to determine a stopping point, since this can
always be calculated readily.

 \section*{Acknowledgments}
 The authors gratefully acknowledge support for aspects of this work under ONR NICOP grant N62909-19-1-2128.

\if\appendixa1
\section*{Appendix A: Projections and eigenvectors}

 For an operator $\B$ with a unique maximum-modulus eigenvalue
 $\lambda$ and corresponding eigenvector $\v$, the pair $(\lambda,\v)$
 can be identified by iteration, by successively applying $\B$ for any
 initial vector $\u $:  
 
 \medbreak
 For large $n$,
 %\rightarrow \infty$, 
 $\B^n\u /||\B^n\u || \sim \v$ and 
  $<\B^{n+1}\u ,\B^n\u >/||\B^n\u || \sim \lambda$  (provided $<\u
 ,\v> \neq 0$) where $<.>$ denotes the inner product, 

\medbreak

If $\B$ were self-adjoint then the component of $\u $ with respect to
$\v$ is just  
$ <\u ,\v>\v \equiv \P(\u )$, say, where 
$\P$ is the rank 1 orthogonal projection onto the line (1-d subspace)
$\mathbb{C}\v$. 
This component can be found, and subtracted from the problem, without
requiring knowledge of 
any other eigenvectors. 
  If $\Q$ is the projection of $\v$ onto the space spanned by the remaining
  eigenvectors - in this case the orthogonal complement of $\v$ - then
  in this case $\Q = 1-\P$, and the resulting vector   
  $\Q\u  = \u  - \P(\u)$ is orthogonal to $\v$.\footnote{
Note that for a non-zero projection $Q$ in Hilbert space the following
are equivalent: 
 	$Q$ orthogonal $\iff Q$ is self-adjoint $\iff ||\Q||=1 \iff
||1-Q||=1 \iff$ range and null-space of $Q$ are orthogonal, etc.} 
 
However, our operator $\B = RL^{-1}$ is not self-adjoint and has
non-orthogonal eigenvectors; thus the 
projection of $\v$ onto the residual eigenspace will not in general be
orthogonal and 
cannot be found simply from $\P$, independently of the remaining
eigenvectors.    

In order to identify and subtract the relevant component:
We first define $C = L^{-1}\B \equiv L^{-1}RL^{-1}$ and compute the
eigenvector matrix $V$ of $C$ together with its eigenvalues
$\{\lambda_i\}$ ordered by decreasing modulus, $|\lambda_1| >
|\lambda_2| > \dots > |\lambda_n|$ say. 
So the columns $\v_k$, say, of $V$ are the eigenvectors of $C$ and the
dominant eigenvector is $\v_1$.  The projection $\Q$ of $\v_1$ onto
the residual eigenspace is not self-adjoint.  

We now calculate the inverse $V^{-1}$ of the eigenvector matrix and
apply it to the incident field $\psinc$, to get $\psi'=V^{-1}\psinc$.
Then $\psinc = \sum_k \psi'_k \v_k$, and the modified incident field
with the $v_1$ component removed is simply $\psinc' = \psinc - \psi'_1
\v_1$.  
More generally, if $\lambda_i$ for $i \le K$ are the dominant
eigenvalues then we can subtract the corresponding eigencomponents to
form the modified field 
\beq \psinc' = \sum_{k> K} \psi'_k \v_k     \label{modified}\eeq

This procedure is carried out only in order to analyse the behaviour
of the system; it is clearly not a suitable method for solving the
original problem since the computational cost is similar to that of
full system.   
\fi

\if\appendixb1
\section*{Appendix B: Higher-dimensional vector Shanks method}

We consider an extension of the principles in section
\ref{vshanks} above  to a larger number of exponentially growing
terms. 
Suppose then that %$\v$, $\w$ are vectors,  and that
\beq S_k = S +\lambda^k \v  \label{v11}  + \delta^k\w ~~~~{\rm for~}k
=1,\dots ~~~~\eeq
where $\lambda \ge \delta$, and suppose that $\lambda$ and $\v$ have
already been identified from the first application of Shanks (or from the
power method).

Again define
\beq \Delta S_k ~=~ S_{k+1} -S_k ~~~~{\rm for~}k =1,\dots ~~~~~.\eeq

Consecutive terms $\Delta\Sn$ are now coplanar rather than colinear, and lie in
the plane generated by $\v,\w$. 
Denote 
\beq D_k  ~=~ S_k - \lambda^k \v  ~\equiv~ S  + \delta^k\w ~~~~{\rm for~}k
=1,\dots ~~~~.   \label{v12} \eeq  
and define $\Delta D_k = D_{k+1}-D_k$. We then obtain 

\beq \Delta D_k =  \Delta S_k -(\lambda^{k+1}-\lambda^k)\v \eeq
and we can write
\beq \Delta D_n ~=~ (\delta^{n+1}-\delta^n)\w  \eeq
\beq \Delta D_{n+1} ~=~ (\delta^{n+2}-\delta^{n+1})\w  
~~=~~ \delta \Delta D_n ,\eeq
so applying the same reasoning as for equations (\ref{v1})-(\ref{v4}) we obtain:
\beq 
\delta = \frac{ \left<\Delta D_{n+1}, \Delta D_n \right>}
{ ||\Delta D_n||^2}. ~~~~~~ \w = \frac{ \Delta D_n} {\delta^{n+1}-\delta^n}   .
\label{v44}     \eeq
This leads to a modified higher-order Shanks transformation with
$S_k,\lambda, \v$ replaced by $D_k,\delta, \w$  (where the dependence
on $\v, \w$ is implicit): 
\beq 
G'_n = \frac{\delta_n D_n - D_{n+1}}{\delta_n-1} \label{vectorshanks2}
\eeq 
\fi

\bigbreak

\if\appendixc1
\section*{Appendix C: Dilating eigenvector solutions}

\def\EQN{\eqno}
\def\beq{\begin{equation}}
\def\eeq{\end{equation}}
\def\beqn{\begin{eqnarray}}
\def\eeqn{\end{eqnarray}}
\def\bef{\begin{figure}}
	\def\eef{\end{figure}}
\def\del{\partial}
\def\Psiz{\Phi}
\def\psinc{\psi_{inc}}
\def\vnc{\v_{inc}}
\def\r{{\bf r}}
\def\rr{{\bf r'}}
\def\d{{\bf d}}
\def\h{{\bf h}}
\def\A{{\bf A}}
\def\B{{B}}   %{{\bf B}}
\def\bea{\begin{aligned}}
	\def\eea{\end{aligned}}
\def\zz{Z}
\def\Arcsin{{\rm Arcsin}}
\def\Arctan{{\rm Arctan}}
\def\t{\theta_n}
\def\tt{\frac{\theta_n}{2}}
\def\hmax{h_{\max}}
\def\hmin{h_{\min}}
\def\P{{P}}   %{{\bf P}}
\def\Q{{Q}}     %{{\bf Q}}
\def\Ecal{{\cal{E}}}
\def\Fcal{{\cal{F}}}
\def\Gcal{{\cal{G}}}
\def\Scal{{\cal{S}}}
\def\Hcal{{\cal{H}}}
\def\Bcal{{\cal{B}}}
\def\Pcal{{\cal{P}}}
\def\Qcal{{\cal{Q}}}
\def\Ical{{\cal{I}}}
\def\Jcal{{\cal{J}}}
\def\Lcal{{\cal{L}}}
\def\Rcal{{\cal{R}}}
\def\Mcal{{\cal{M}}}
\def\Ncal{{\cal{N}}}
\def\Einc{{\u}}
\def\Escat{{E_s}}
\def\Epartial{\frac{\partial E}{ \partial n}}

\def\Sn{{S_n}}
\def\Sone{{S_{n+1}}}
\def\Stwo{{S_{n+2}}}
\def\v{{\bf v}}
\def\u{{\bf u}}
\def\w{{\bf w}}

\def\rbar{{\bf r}}
\def\lbar{{\bar{\lambda}}}
\def\J{{\bf J}}
\def\t{{\bf t}}

\parindent 0 true cm

\font\eightrm=cmr8
\subsection*{Alternative text for section \ref{projectionmethod}:}

\medbreak
Suppose $A$ is a [possibly invertible] linear operator on a Hilbert space
and that we wish to solve $A\u = \v$  where 
$\v$ is a given vector,  to find $\u ~=~ A^{-1}\v$.  
Let $A=L+R$ be any decomposition into bounded linear operators $L$,
$R$  where $L$ is invertible.  

Where this solution exists it can formally be written as a series:
\beq
A^{-1}\v ~=~ L^{-1} \sum_{n=0}^\infty B^n \v
\label{(3.6new)} 
\eeq
where $\B = - R L^{-1}$.
We also consider the partial sums of (finite approximations to) this  
series:
\beq S_n (\v)=  L^{-1}\sum_{k=0}^n B^k \v
\label{partial}
\eeq 

Of course, this need not converge even when $A^{-1}\v$
is well-defined.   In particular 
$B$ may  have eigenvectors whose eigenvalues have moduli greater than
1 (we'll call these `dilating').
Any components of $\v$ in such directions will cause the RHS of
(\ref{(3.6new)}) 
to diverge.

\medbreak
We examine here how the exact
solution relates to the series solution for such components.
Specifically we ask: (1) What
is the exact solution $A^{-1}\v_0$ at a dilating
eigenvector $\v_0$ of $B$?    
(2) Can we obtain this solution from the $n$-th iterates or the
partial sums?  

\medbreak
---------------------------------
\medbreak
Suppose then that $(\v_0,\lambda_0)$ is an eigenvector/eigenvalue pair
for $B$.  Provided 
$\lambda_0 \neq 1$, expression (\ref{partial}) then reduces to   
\beq
L^{-1} \sum_{k=0}^n B^k \v_0
~=~ L^{-1} \sum_{k=0}^n \lambda_0^k \v_0
~=~  \left(\frac{\lambda_0^{n+1}-1}{\lambda_0-1}\right)~L^{-1}~\v_0
\label{finitesum}
\eeq
which clearly converges for $|\lambda_0| < 1$ and diverges for $|\lambda_0|>1$.
We can show that the correct limit in {\it both} cases is obtained by
discarding the term $\lambda_0^{n+1}$ in the numerator on the right:

\medbreak
Since $B=-RL^{-1}$ and $B\v_0 = \lambda_0\v_0$ we have
\beq (1-\lambda_0)\v_0 ~~=~~\v_0 + RL^{-1}\v_0 ~~=~~ (L+R)L^{-1}\v_0
~~=~~ AL^{-1} \v_0 
\label{lemma1a}\eeq 
so that, as claimed,
\beq A^{-1}\v_0 ~~=~~ \frac{1}{1-\lambda_0}L^{-1}\v_0  . \label{lemma1b}\eeq

\medbreak
(Comparing (\ref{lemma1b}) and (\ref{finitesum}) we see that for $|\lambda_0| < 1$ the
series converges automatically to the correct limit,
but in for dilating eigenvectors the term $\lambda_0^{n+1}$ in the
numerator must be discarded explicitly.)

\bigbreak

It is straightforward to extend this to multiple eigenvectors:
Suppose now that we have a [not necessarily complete] set of eigenvectors 
$\v_i$ for $B = -RL^{-1}$ with  
eigenvalues $\lambda_i$, and that these eigenvectors have 
multiplicity one.   

Denote by $P_i$ the rank 1 projection  onto the eigenspace of $\v_i$, so that
$P_iP_j = 0$ for $i\neq j$ and 
the sum of any subset of $\{P_i\}$ is a projection.  

{\sl [Remark: If $B$ were self-adjoint these eigenvectors would form an
	orthogonal set.  
	However, this is not the case in our application. Thus the range
	and null space of any $P_i$ are not in general orthogonal subspaces
	\footnote{
		Note that for a nonzero projection $Q$ in Hilbert space the following
		are equivalent: 
		$Q$ orthogonal $\iff Q$ is self-adjoint $\iff ||\Q||=1 \iff
		||1-Q||=1 \iff$ range and null-space of $Q$ are orthogonal, etc.}.]}

For a general vector $\v$ in the span of the eigenvectors, the partial sums 
$S_n (\v)$ 
can be written in terms of the projections $P_i$ as

\beq S_n (\v)~=~ L^{-1}\sum_i \sum_{k=0}^n \left(\lambda_i^k P_i \right)\v 
~=~  L^{-1}\sum_i\sum_{k=0}^n \left(\lambda_i^k \alpha_i\v_i \right)\eeq 
where $P_i\v = \alpha_i\v_i$, which gives
\beq 
S_n (\v)~=~ L^{-1}\sum_i \alpha_i \Big(\frac{ \lambda_i^{n+1} - 1
}{\lambda_i - 1} \Big)~\v_i   . 
\label{snv}
\eeq

So for any such $\v$  
\beq A^{-1}\v ~=~
\sum_i \Big(\frac{ 1}{1-\lambda_i} \Big)~L^{-1} P_i \v_i. 
\label{exact}\eeq

\fi

\bibliographystyle{unsrt}

\bibliography{refs}

\end{document}